\documentclass[a4paper,11pt]{article}
\usepackage[english]{babel}
\usepackage{amsmath}
\usepackage{siunitx}
\usepackage[a4paper,left=3cm,right=3cm,top=2.5cm,bottom=2.5cm]{geometry}
\usepackage{hyperref}
\hypersetup{colorlinks=true, linkcolor=blue, urlcolor=red,
pdftitle={hyper_roms_for_moving_domains}}
\usepackage{fancyvrb}
\usepackage{fancyhdr, lastpage}
\usepackage{bm}
\usepackage{pdfpages}
\usepackage{authblk}
\usepackage{graphicx}
\usepackage{float}
\usepackage{comment}
\usepackage{caption}
\usepackage{cite}
\usepackage[colorinlistoftodos]{todonotes}
\usepackage{mathtools}
\usepackage{amssymb}
\usepackage{amsthm}
\usepackage{multirow}
\usepackage{tabularx}
\usepackage{longtable}
\usepackage{setspace}
\usepackage{xcolor}
\usepackage{titlesec}
\usepackage{booktabs}
\usepackage[draft,inline,marginclue]{fixme}
\FXRegisterAuthor{gs}{ags}{\color{red}GS}

\DeclareGraphicsExtensions{.pdf,.png}
\usepackage{algorithm, algorithmicx, algpseudocode}
\usepackage{cleveref}
\fxsetup{status=draft}
\graphicspath{ {./images/} }

\usepackage{todonotes}
\DeclareRobustCommand{\intodo}[1]{%
  \todo[inline]{#1}%
}

\addtocontents{toc}{\begingroup\protect\renewcommand{\protect\intodo}[1]{}}
\AtEndDocument{%
  \addtocontents{toc}{\endgroup}
}

\newcommand{\RA}{\color{blue}}
\newcommand{\RB}{\color{red}}

\title{An efficient hyper-reduced-order model for segregated solvers for geometrical parametrization problems}
\author[1]{Valentin Nkana Ngan\footnote{vkanang@sissa.it}}
\author[2]{Giovanni Stabile\footnote{giovanni.stabile@santannapisa.it}}
\author[3]{Andrea Mola\footnote{andrea.mola@imtlucca.it}}
\author[1]{Gianluigi Rozza\footnote{grozza@sissa.it}}

\affil[1]{Mathematics Area, mathLab, SISSA, via Bonomea 265, I-34136 Trieste, Italy}
\affil[2]{The Biorobotics Institute, Sant'Anna School of Advanced Studies, V.le R. Piaggio 34, 56025, Pontedera, Pisa - Italy}
\affil[3]{MUSAM Continuum Mechanics Laboratory, Scuola IMT Alti Studi Lucca - Piazza S. Ponziano, 6 - 55100 Lucca, LU, Italy }
\date{\today} 

\begin{document}
\tableofcontents
\maketitle
\begin{abstract}
We propose an efficient hyper-reduced order model (HROM) designed for segregated finite-volume solvers in geometrically parametrized problems. 
The method follows a discretize-then-project strategy: the full-order operators are first assembled using finite volume or finite element discretizations and then projected onto low-dimensional spaces using a small set of spatial sampling points, selected through hyper-reduction techniques such as DEIM. 
This approach removes the dependence of the online computational cost on the full mesh size. 
The method is assessed on three benchmark problems: a linear transport equation, a nonlinear Burgers’ equation, and the incompressible Navier--Stokes equations. 
The results show that the hyper-reduced models closely match full-order solutions while achieving substantial reductions in computational time. 
Since only a sparse subset of mesh cells is evaluated during the online phase, the method is naturally parallelizable and scalable to very large meshes. 
These findings demonstrate that hyper-reduction can be effectively combined with segregated solvers and geometric parametrization to enable fast and accurate CFD simulations.
\end{abstract}

\listoffixmes
\section{Motivation and state-of-the-art}

Many engineering design and optimization problems require the repeated solution of partial differential equations (PDEs) for different parameter values. 
These parameters include physical quantities, such as boundary conditions, engineering properties, such as material coefficients, and geometric variations. 
Geometric parametrization is essential in applications including biomedical devices \cite{ballarin2016} and aerodynamic shape design \cite{quarteroni2015reduced}, but it significantly increases computational cost due to mesh deformation or regeneration. 
The large number of parameters, their nonlinear interactions, and their potentially high dimensionality make such simulations challenging.

Model order reduction techniques have been developed to alleviate these costs. 
Reduced basis methods (RBMs) approximate the solution manifold of parametrized PDEs using low-dimensional spaces constructed from high-fidelity solution snapshots. 
RBMs have been successfully applied to a wide range of problems, including linear elasticity and material science \cite{patera2007reduced}, heat and mass transfer \cite{grepl2007efficient}, acoustics \cite{buchan2013reduced}, potential flows \cite{rozza2005reduced}, and microfluid dynamics \cite{lassila2014model}. 
For complex multiphysics systems, the most effective RBM formulations rely on projection-based reduced-order models (ROMs), which systematically project the governing equations onto reduced spaces.

Despite their success, projection-based ROMs often fail to deliver substantial computational speed-ups in realistic engineering settings. 
This limitation is particularly evident in problems involving geometric parametrization, moving domains, or strong nonlinearities, where the assembly and evaluation of reduced operators still depend on the full-order discretization size \cite{zucatti2021calibration}. 
As a result, the online cost may remain dominated by full-order operations, limiting the practical impact of the reduction.

Hyper-reduced order methods (HROMs) were introduced to address this issue by further reducing the cost of operator evaluation. 
Classical approaches, such as empirical interpolation and its variants, approximate nonlinear operators using a small number of spatial sampling points or quadrature evaluations \cite{ryckelynck2005priori}. 
While these methods have proven effective in many finite element settings, several limitations remain for their use in large-scale finite-volume simulations with geometric variability.

First, most existing HROM formulations are developed for monolithic solvers and assume a tight coupling between all governing equations. 
This assumption is not compatible with segregated solution strategies, such as SIMPLE, PISO, or PIMPLE, which are widely used in industrial CFD codes. 
Second, hyper-reduction techniques are often presented independently of the solver structure, making their integration into existing finite-volume frameworks non-trivial. 
Third, in geometrically parametrized problems, the dependence of discrete operators on mesh topology and connectivity is rarely addressed explicitly, limiting robustness under mesh motion.

The present work addresses these gaps by introducing a hyper-reduced order modeling strategy specifically designed for segregated finite-volume solvers and geometrically parametrized domains. 
The proposed method follows a discretize-then-project approach, in which the full-order finite-volume operators are assembled first and then evaluated online using a small, carefully selected subset of mesh cells. 
This design removes the dependence of the online computational cost on the full mesh size, while preserving the structure of segregated solvers and accommodating geometric variations.

The remainder of the paper is organized as follows. 
The projection-based reduced-order formulation is introduced in \autoref{sec:galproj}, followed by the proposed hyper-reduction strategy in \autoref{sec:hyper_roms}. 
Numerical results and performance analyses are presented in \autoref{resdiscuss} for three benchmark problems: a two-dimensional transport equation, a two-dimensional Burgers equation, and a two-dimensional laminar flow around a cylinder. 
Conclusions and perspectives are given in \autoref{conclusion}.

\section{Hyper-reduction}

\subsection{Galerkin projection}
\label{sec:galproj}

We first recall the standard Galerkin projection used to construct projection-based reduced-order models. 
The presentation is kept general and applies to different fields of interest, such as velocity, pressure, or temperature.

We consider a linear system obtained from the spatial discretization of a system of PDEs, for instance the momentum, energy, or pressure Poisson equations,
\begin{equation}
\label{linSystem}
\mathbf{A}(\mu)\,\bm{w}(\bm{x},t;\mu) = \bm{b}(\mu),
\end{equation}
where $\mathbf{A}(\mu)\in\mathbb{R}^{N_h\times N_h}$ is the discrete operator, $\bm{b}(\mu)\in\mathbb{R}^{N_h}$ is the source term, and $\bm{w}\in\mathbb{R}^{N_h}$ denotes the discrete solution field.

In projection-based reduced-order modelling, the solution is approximated in a low-dimensional subspace spanned by $N_r \ll N_h$ basis functions,
\begin{equation}
\label{eq:separability}
\bm{w}(\bm{x},\cdot)\simeq \sum_{i=1}^{N_r} a_i(\cdot)\,\bm{\phi}_i(\bm{x})
= \bm{\Phi}\bm{a},
\end{equation}
where $\bm{\phi}_i$ are typically obtained from Proper Orthogonal Decomposition, collected in the matrix $\bm{\Phi}\in\mathbb{R}^{N_h\times N_r}$, and $\bm{a}\in\mathbb{R}^{N_r}$ contains the reduced coefficients.

Substituting this approximation into \autoref{linSystem} and enforcing Galerkin orthogonality yields the reduced system
\begin{equation}
\label{redSystem}
\mathbf{A}^r \bm{a} = \bm{b}^r,
\end{equation}
with $\mathbf{A}^r = \bm{\Phi}^\intercal \mathbf{A}\bm{\Phi}\in\mathbb{R}^{N_r\times N_r}$ and $\bm{b}^r = \bm{\Phi}^\intercal \bm{b}\in\mathbb{R}^{N_r}$.
The resulting system is dense but low dimensional and can be solved efficiently using standard linear algebra techniques, such as rank-revealing QR factorizations \cite{eigenweb}.

While this procedure reduces the number of unknowns, the construction of $\mathbf{A}^r$ and $\bm{b}^r$ still requires operations that scale with the full-order dimension $N_h$, which limits the achievable speed-up in practice.

\subsection{Hyper-reduction algorithm for geometrical parametrization problems}
\label{sec:hyper_roms}

To remove the dependence of the online computational cost on the full-order dimension $N_h$, we introduce a hyper-reduction strategy tailored to segregated finite-volume solvers, such as SIMPLE, PISO, and PIMPLE.
The goal is to evaluate the reduced operators using only a small subset of mesh cells, while preserving the structure of the underlying solver.

The key idea is to approximate the action of the operator $\mathbf{A}(\mu)\bm{\Phi}$ by evaluating it only at a selected set of spatial locations, often referred to as magic or optimal cells.
We assume that a set of sampling indices
\[
\mathcal{J} = \{j_1,\dots,j_s\}\subset\{1,\dots,N_h\}, \qquad s \ll N_h,
\]
is available, obtained for instance using DEIM-type techniques
\cite{aanonsen2009empirical,Drma2016,chaturantabut2010nonlinear,nguyen2008best,chen2021eim,Zimmermann2024,chen2019l1,Dimitriu2017,Lauzon2024,Zimmermann2016}.
The associated selector matrix is defined as $\bm{\mathcal{P}}=[\bm{e}_{j_1},\dots,\bm{e}_{j_s}]\in\mathbb{R}^{N_h\times s}$.

Instead of solving the full Galerkin system, the reduced coefficients are obtained by minimizing the residual of the projected equations in a weighted discrete $L^2$ norm,
\begin{equation}
\label{overdeterminated}
\min_{\bm{a}\in\mathbb{R}^{N_r}}
\left\| \mathbf{A}(\mu)\bm{\Phi}\bm{a}(\mu) - \bm{b}(\mu) \right\|_{\bm{S}}^2,
\end{equation}
where $\bm{S}$ is a diagonal matrix containing weights associated with the spatial discretization.
In the finite-volume setting, these weights correspond to the inverse cell volumes, $S_{ii}=1/\lvert\Omega_i\rvert$.

Direct evaluation of \autoref{overdeterminated} remains expensive, since it involves the full product $\mathbf{A}(\mu)\bm{\Phi}$.
To avoid this cost, we restrict the residual to the sampled cells and reformulate the problem as a masked, weighted least-squares system,
\begin{equation}
\label{leastsquare}
\min_{\bm{a}\in\mathbb{R}^{N_r}}
\left\|
\bm{\mathcal{P}}^\intercal \sqrt{\bm{S}}\,\mathbf{A}(\mu)\bm{\Phi}\bm{a}(\mu)
-
\bm{\mathcal{P}}^\intercal \sqrt{\bm{S}}\,\bm{b}(\mu)
\right\|^2.
\end{equation}

The reduced operator is assembled by extracting only the rows of $\mathbf{A}(\mu)$ corresponding to the sampling indices and exploiting the sparsity of the finite-volume stencil.
As a result, the computational cost scales as $\mathcal{O}(s\,r\,N_r)$, where $r$ is the average number of non-zero entries per row of $\mathbf{A}(\mu)$, and is completely independent of $N_h$.
Once the reduced system is constructed, all remaining operations depend only on the reduced dimensions.

This procedure preserves the segregated structure of the original solver and naturally accommodates geometric parametrization, since only local stencil information around the sampled cells is required. The full algorithm is summarized in Algorithm~\ref{alg:hyper_reduced_algo}.

\begin{algorithm}
\caption{Hyper-reduced order algorithm}
\label{alg:hyper_reduced_algo}
\begin{algorithmic}[1]
\Require{$\bm{A} \in \mathbb{R}^{N_h \times N_h}$ (row-major), 
  $\bm{b} \in \mathbb{R}^{N_h}$, $\mathbf{V}\in \mathbb{R}^{N_h}$,
  $\mathcal{J}$, $\mathbf{\Phi} \in \mathbb{R}^{N_h \times N_r}$;}
\Ensure{$\mathbf{A}^r \in \mathbb{R}^{s \times N_r}$, $\bm{b}^r \in \mathbb{R}^{s}$}\; \Comment{$s \geq N_r$}
\State $s \gets |\mathcal{J}|$\; \Comment{number of magic points}
\State $N_r \gets$ number of columns of $\mathbf{\Phi}$\; \Comment{number of modes}
\For{$k \gets 1$ to $s$}
  \State $i \gets \mathcal{J}[k]$ \Comment{row index from magic points}
  \State $\bm{z} \gets \bm{0}\in \mathbb{R}^{N_r}$ \;
  \For{non-zero entry $(i,j,\bm{A}_{ij})$ in row $i$ of $\bm{A}$}
     \State $\bm{z} \gets \bm{z} + \bm{A}_{ij} \cdot \mathbf{\Phi}[:,j]$ \;
  \EndFor
  \State $v \gets \frac{1}{\mathbf{V}[i]}$\; \Comment{reciprocal of mesh volume}
  \State $\mathbf{A}^r[k,:] \gets \sqrt{v}\,\bm{z}$ \;
  \State $b^r[k] \gets \sqrt{v}\, \bm{b}[i]$ \;
\EndFor
\end{algorithmic}
\end{algorithm}

\section{Case studies,  simulation results, and discussion}
\label{resdiscuss}
In engineering applications, the matrix $\mathbf{A}$ that defines the least-squares problem in Proposition 2.2
may depend on a parameter vector $\mu$ (e.g., parameters that account for varying system properties such as material constants or boundary conditions).
	
\subsection{Basic performance tests}
\begin{figure}
\centering
\includegraphics[width=0.95\textwidth]{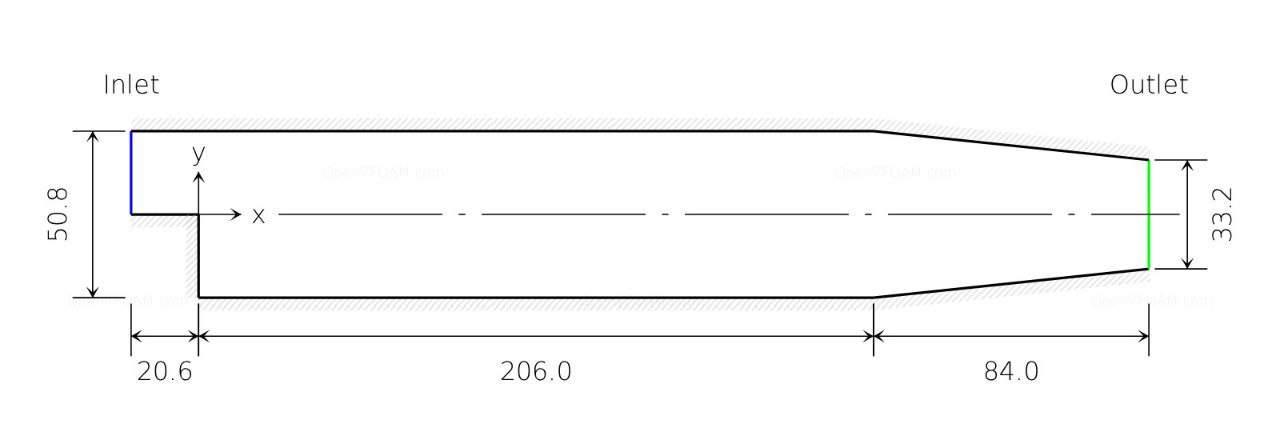}
	\caption{The geometry of the domain}	
   \label{Backward_geo}
\end{figure}

\begin{figure}
\centering
\includegraphics[width=0.95\textwidth]{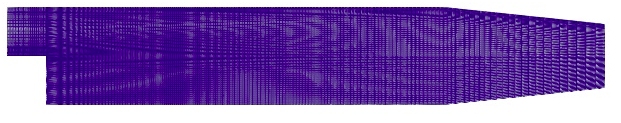}
	\caption{The mesh domain}	
   \label{mesh_geo}
\end{figure}

For simple performance tests, we present results on the 2D scalar transport equation and 2D Burgers equation.
The geometry used in these test cases is the backward-facing step, as shown in \autoref{Backward_geo}. 
We solve the 2D scalar transport equation and 2D Burgers' equation with OpenFOAM \cite{OpenFOAM}. 
We employ the FVM in a structured orthogonal grid of 12225 cells as shown in \autoref {mesh_geo}.  
At the inlet boundary, non-homogeneous Dirichlet and zero gradient conditions are prescribed for the temperature or the velocity field. 
At the outlet boundary, zero gradient is prescribed for $T$ or $\boldsymbol{u}$.
On the upper and lower walls (top and bottom) no no-slip boundary conditions are prescribed for the temperature.  
In the full-order simulation, the Gauss linear scheme was selected for the approximation of the
gradients, and the Gauss linear scheme with \textit{non-orthogonal correction} was selected to approximate the Laplacian term. 
A Gauss linear scheme was instead used for the approximation of the convective term. 
The FOM numerical scheme is the SIMPLE algorithm. The time discretization is performed with the semi-implicit Euler method.   
The linear solver used for the resulting linear system is based on the GAMG with Gauss-Seidel smoother until a tolerance of 1e-08 on the finite volume method residual is reached.

\subsubsection{2D scalar transport equation}
This test case is designed to assess the effectiveness of the proposed
hyper-reduction strategy for linear convection–diffusion operators discretized
with a segregated finite-volume scheme. In particular, it evaluates whether
accurate solutions can be obtained while reducing the online computational cost
to a level independent of the full mesh size.

We consider the two-dimensional transport equation \autoref{eq:transport} on the
backward-facing step geometry shown in \autoref{Backward_geo}. The problem is
discretized using the finite volume method on a structured orthogonal mesh
consisting of 12\,225 cells, as depicted in \autoref{mesh_geo}. Time integration
is performed using a semi-implicit Euler scheme with a time step
$\Delta t = 10^{-4}\,\mathrm{s}$ up to a final time $t_f = 0.25\,\mathrm{s}$.
The kinematic viscosity is set to $\nu = 4 \times 10^{-5}\,\mathrm{m}^2/\mathrm{s}$,
and snapshots are collected at every time step.
\begin{equation}
\label{eq:transport}
\begin{cases}
\frac{\partial T }{\partial t}  +   \nabla (\boldsymbol{u}T) - \nu \Delta T  = S_{T}, \quad (\bm{x}, t) \in \Omega \times [0, T_f] \\
T(\bm{x}, 0) = T_{0}(\bm{x}), \quad \bm{x} \in \Omega \\ 
T(\bm{x}, t)|_{\partial \Omega} = 0,   \quad  t \in [0, T_f].
\end{cases}
\end{equation}
Where $\Omega \subset \mathbb{R}^2$ is a bounded domain with boundary  $\partial \Omega$.  
$T$ being the temperature, $S_{T}$ source term here set to zero.  $\nu$ the viscosity, and $\boldsymbol{u}$ velocity. 

\begin{figure}
    \centering
    \includegraphics[width=0.95\linewidth]{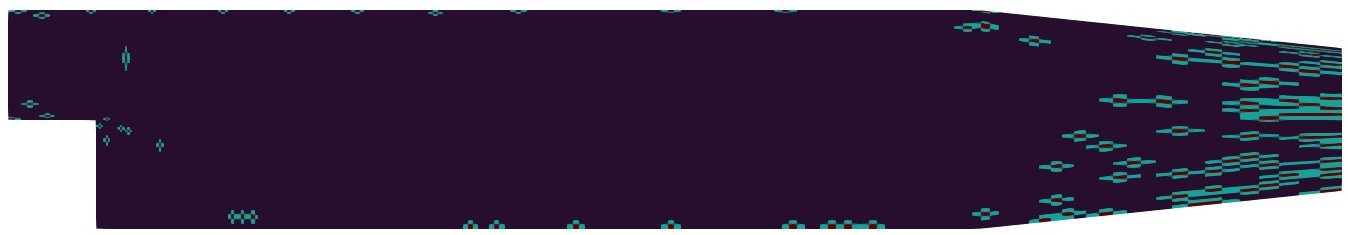}
    \caption{The 100 DEIM points are represented in red, the stencils of the cells and associated degrees of freedom needed for the evaluation of the discrete differential operators are in light-blue. 
    The discarded nodes in the evolution of the dynamics are in blue. The stencil is made of 1 layer of cells}
    \label{fig:transport_optim_points}
\end{figure}

The hyper-reduction is performed using a DEIM-based sampling strategy applied to
the reduced operators. \autoref{fig:transport_optim_points} illustrates the
selected interpolation points together with the associated finite-volume
stencils required for evaluating the discrete differential operators. Only a
single layer of neighbouring cells is retained, highlighting the strictly local
nature of the proposed stencil-based evaluation.

Table~\ref{tab:OffOnlineTransport} reports the offline and online computational
times for different numbers of temperature modes. For $N_T = 7$ modes, the
hyper-reduced model achieves a speed-up of approximately $7.3\times$ compared
to the full-order simulation, while maintaining similar accuracy for $N_T = 5$
modes. These results confirm that, for linear problems, the online cost of the
reduced model is dominated by operations on the sampled cells and is therefore
independent of the total number of mesh elements.

\begin{table}
\centering
\begin{tabular}{|p{2cm}||p{5cm}||p{3cm}||p{2cm}|}
	\hline
	\textbf{Stages} & \textit{\# of modes} & Time [s] & {\RA Speed-up $[T_{Off}/ T_{On}]$} \\
 \hline
	\textbf{Offline} &  - & 109.211 & -\\
 \hline
    \multirow{3}{*}{\textbf{Online}} & $N_T=7$ & 14.9101 & {7.324}\\
       \cline{2-4}
       &    $N_T= 5 $ & 13.7344 & { 7.951}\\
       \hline
\end{tabular}
\caption{\RB Speed-up, offline and online times comparison varying the number of modes.}
\label{tab:OffOnlineTransport}
\end{table}

A qualitative comparison between the full-order and hyper-reduced solutions is
shown in \autoref{fig:T_comparison}. The temperature fields produced by the
hyper-reduced model closely match the reference solutions at all reported time
instances, and the associated error fields remain localized and small. This
indicates that the reduced basis, combined with sparse operator evaluation, is
sufficient to capture the dominant transport dynamics.

\begingroup
\setlength{\tabcolsep}{0.05pt} 
\begin{figure}
\centering
\begin{tabular}{ccc}
\includegraphics[width=0.33\linewidth]{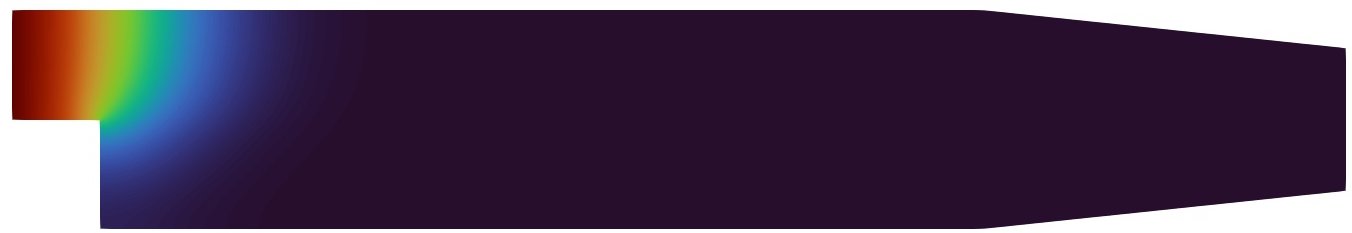} &  \includegraphics[width=0.33\linewidth]{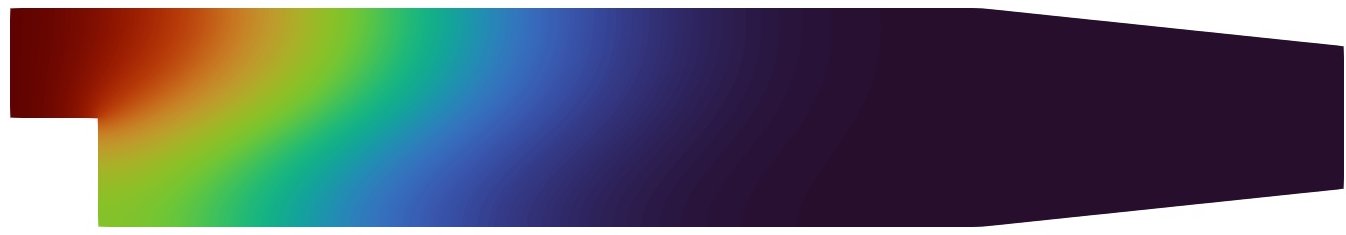} &
\includegraphics[width=0.33\linewidth]{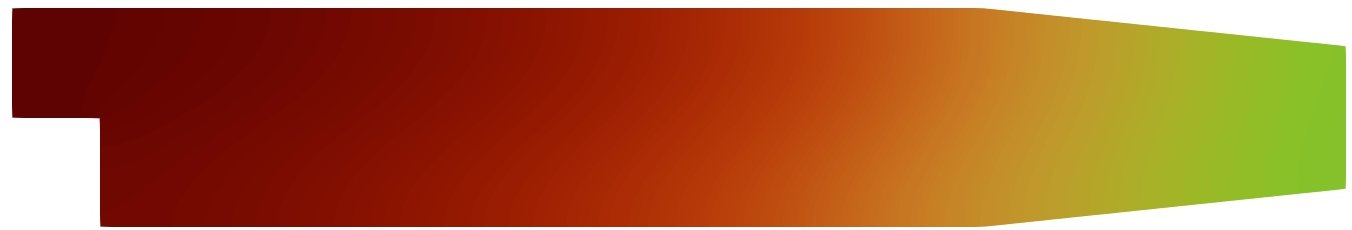}\\
\includegraphics[width=0.23\linewidth]{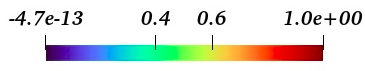}  &  \includegraphics[width=0.23\linewidth]{images/T_scale}  & \includegraphics[width=0.23\linewidth]{images/T_scale} \\
\includegraphics[width=0.33\linewidth]{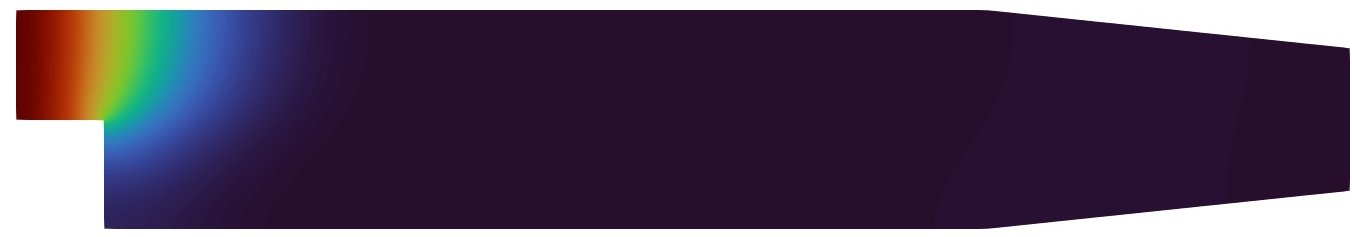} &  \includegraphics[width=0.33\linewidth]{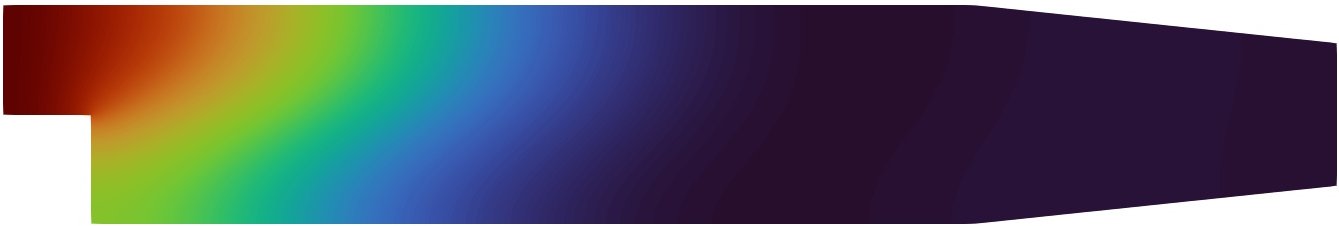} &
\includegraphics[width=0.33\linewidth]{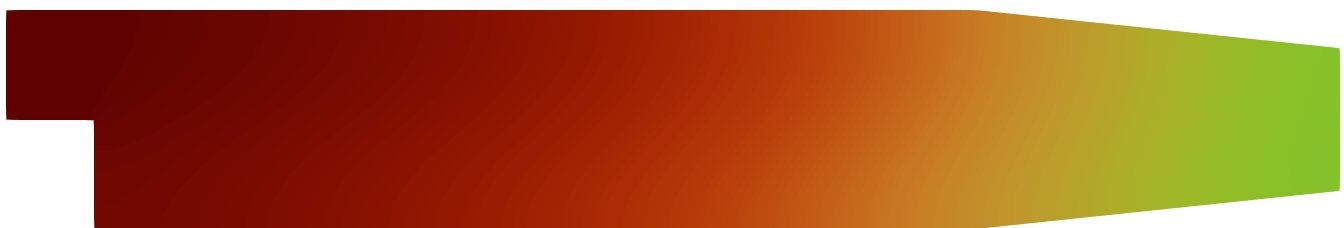}\\
\includegraphics[width=0.23\linewidth]{images/T_scale}  &  \includegraphics[width=0.23\linewidth]{images/T_scale}  & \includegraphics[width=0.23\linewidth]{images/T_scale} \\
\includegraphics[width=0.33\linewidth]{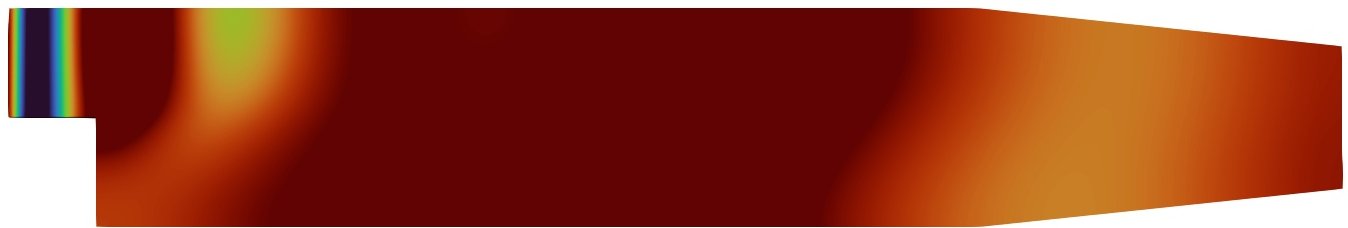} &  \includegraphics[width=0.33\linewidth]{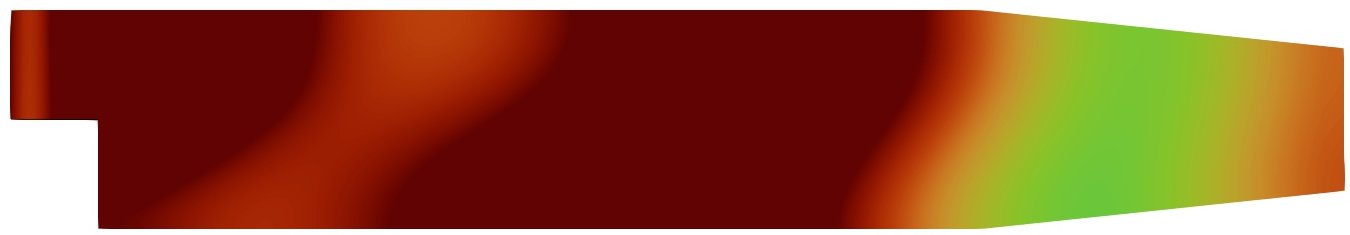} &
\includegraphics[width=0.33\linewidth]{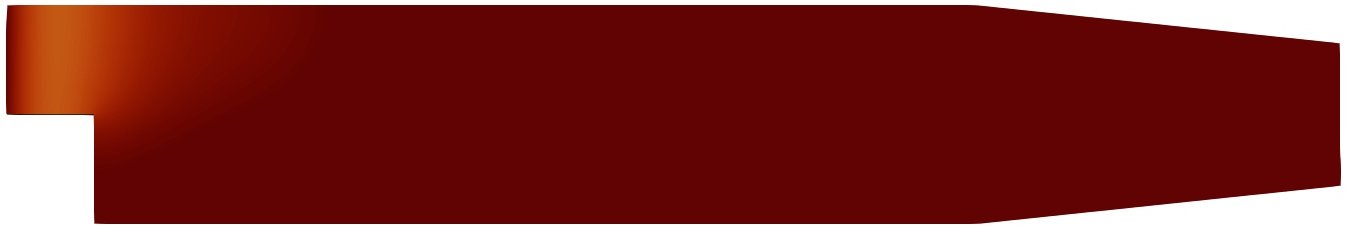}\\
\includegraphics[width=0.23\linewidth]{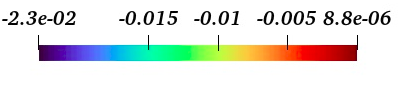}  &  \includegraphics[width=0.23\linewidth]{images/Te_scale}  & \includegraphics[width=0.23\linewidth]{images/Te_scale} \\
\end{tabular}
\caption{Comparison of the temperature field. First row original solutions, second row hyper-reduced solutions using 7 \textbf{optimal points}, and third row the error associated. 
Snapshots are taken at $t=0.002\ $s, $t=0.01\ $s, and $t=0.052\ $s.}	
 \label{fig:T_comparison}
\end{figure}
\endgroup
Overall, this test case demonstrates that the proposed hyper-reduction strategy
can be seamlessly integrated into a segregated finite-volume framework and
deliver substantial computational savings for linear convection–diffusion
problems without compromising solution quality. It also establishes a baseline
against which the more challenging nonlinear and coupled flow problems are
assessed in the following sections.

\subsubsection{2D Burgers Equation}
This test case is introduced to assess the performance and limitations of the
proposed hyper-reduction strategy in the presence of strong nonlinear dynamics.
The two-dimensional Burgers equation represents a prototypical nonlinear
convection–diffusion problem and serves as a simplified model for nonlinear
advection terms appearing in the Navier--Stokes equations.

We consider the two-dimensional Burgers equation given in
\autoref{eq:Burgers} on the backward-facing step geometry shown in
\autoref{Backward_geo}. The problem is discretized using the finite volume method
on the same structured orthogonal mesh employed for the transport equation,
consisting of 12\,225 cells. Time integration is carried out using a
semi-implicit Euler scheme with a time step $\Delta t = 10^{-4}\,\mathrm{s}$ up to a final time $t_f = 0.15\,\mathrm{s}$. 
\begin{align}
\label{eq:Burgers}
\begin{cases}
& \frac{\partial \boldsymbol{u} }{\partial t}  +  (\boldsymbol{u}\cdot \nabla)\boldsymbol{u } - \nu \Delta \boldsymbol{u } = 0, \quad (\bm{x}, t) \in \Omega \times [0, T] \\
& \bm{u}(\bm{x}, 0) = \bm{u}_{0}(\bm{x}), \quad \bm{x} \in \Omega\\
& \bm{u}(\bm{x}, t)|_{\partial \Omega} = 0  \quad  t \in [0, T].
\end{cases}
\end{align}
Where $\Omega \subset \mathbb{R}^2$ is a bounded domain with boundary  $\partial \Omega$.  $\boldsymbol{u}$ being the velocity, $\nu$ the viscosity. 
The kinematic viscosity is set to
$\nu = 4 \times 10^{-5}\,\mathrm{m}^2/\mathrm{s}$, and solution snapshots are collected every $3\,\Delta t$.

\begin{figure}
    \centering
    \includegraphics[width=0.95\linewidth]{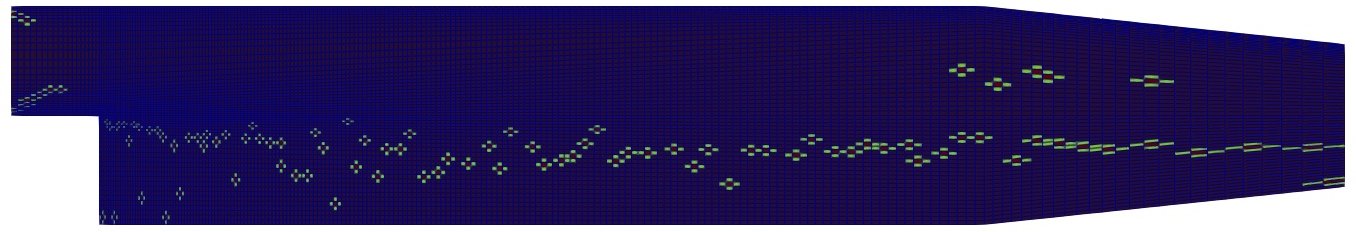}
    \caption{The first 100 DEIM points for the  Burgers equation are represented in red, the stencils of the cells and associated degrees of freedom needed for the evaluation of the discrete differential operators are in light-green. 
    The discarded nodes in the evolution of the dynamics are in blue. The stencil is made of 1 layer of cells}
    \label{fig:burgers_optim_points}
\end{figure}

Table~\ref{tab:OfflineOnlineBurgers} reports the offline and online computational
times for different numbers of velocity modes. The observed speed-up factors,
which remain below $2\times$, are significantly smaller than those obtained for
the linear transport equation. This behavior is expected: for nonlinear
problems, the evaluation of reduced operators still requires repeated assembly
of nonlinear residuals and Jacobian-like contributions at each time step,
thereby limiting the achievable online acceleration.

\begin{table}
\centering
\begin{tabular}{|p{2cm}||p{5cm}||p{3cm}||p{2cm}|}
	\hline
	\textbf{Stages} & \textit{\# of modes} & Time [s] & {\RA Speed-up $[T_{Off}/ T_{On}]$} \\
 \hline
	\textbf{Offline} &  - & 38.7407 & -\\
 \hline
    \multirow{3}{*}{\textbf{Online}} & $N_{\boldsymbol{u}} = 4 $ & 23.1091 & 1.68\\
       \cline{2-4}
       &    $N_{\boldsymbol{u}} = 3 $ & 22.9175& 1.7 \\
       \hline
\end{tabular}
\caption{Speed-up, offline and online times comparison varying the number of modes.}
\label{tab:OfflineOnlineBurgers}
\end{table}

Despite the modest speed-up, the hyper-reduced solutions remain in close
agreement with the full-order reference solutions. \autoref{fig:Ucomparison}
shows that the dominant flow features are accurately captured at all reported
time instances, while the error fields remain bounded and spatially localized.
This indicates that the reduced basis retains sufficient expressive power even
when the dynamics are governed by strong nonlinear interactions.

\begingroup
\setlength{\tabcolsep}{0.05pt} 
\begin{figure}
\centering
\begin{tabular}{ccc}
\includegraphics[width=0.33\linewidth]{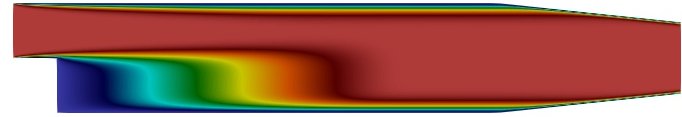} &  \includegraphics[width=0.33\linewidth]{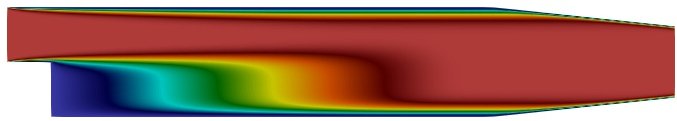} &
\includegraphics[width=0.33\linewidth]{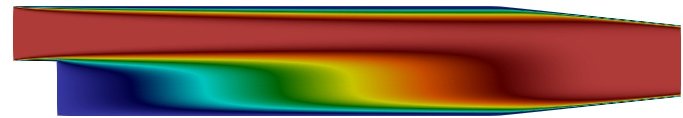}\\
\includegraphics[width=0.23\linewidth]{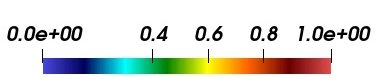}  &  \includegraphics[width=0.23\linewidth]{images/U_fom_scale.jpg}  & \includegraphics[width=0.23\linewidth]{images/U_fom_scale.jpg} \\
\includegraphics[width=0.33\linewidth]{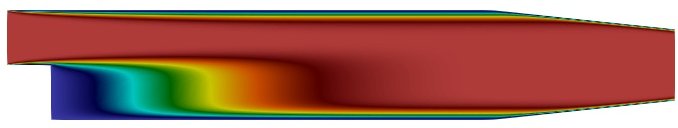} &  \includegraphics[width=0.33\linewidth]{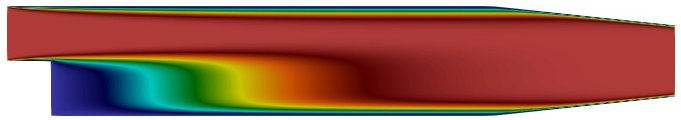} &
\includegraphics[width=0.33\linewidth]{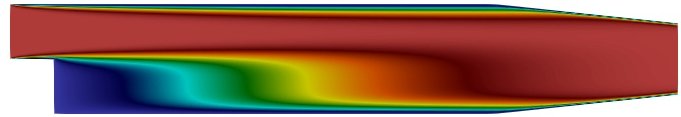}\\
\includegraphics[width=0.23\linewidth]{images/U_fom_scale.jpg}  &  \includegraphics[width=0.23\linewidth]{images/U_fom_scale.jpg}  & \includegraphics[width=0.23\linewidth]{images/U_fom_scale.jpg} \\
\includegraphics[width=0.33\linewidth]{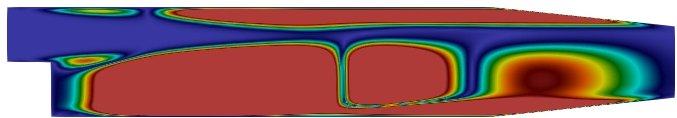} &  \includegraphics[width=0.33\linewidth]{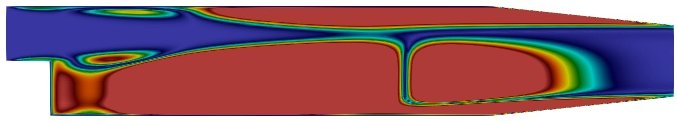} &
\includegraphics[width=0.33\linewidth]{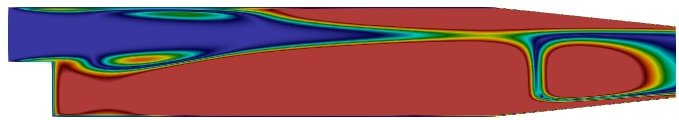}\\
\includegraphics[width=0.23\linewidth]{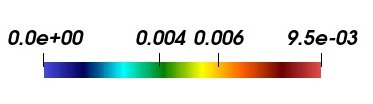}  &  \includegraphics[width=0.23\linewidth]{images/Ue_scale.jpg}  & \includegraphics[width=0.23\linewidth]{images/Ue_scale.jpg} \\
(a) & (b)  &  (c)
\end{tabular}
\caption{Comparison of the velocity field. First row original solutions, second row hyper-reduced solutions using the first \textbf{4 DEIM points}, and third row the associated errors. 
Snapshots are taken at $t=0.0006\ $s, $t=0.03\ $s, and $t=0.156\ $s.}	
 \label{fig:Ucomparison}
\end{figure}
\endgroup
Overall, this test case highlights a fundamental trade-off inherent to
hyper-reduced models for nonlinear systems: while the achievable speed-up is
lower than in linear settings, the proposed approach preserves accuracy and
maintains mesh-independent online complexity. These results justify the use of
the method as a robust building block for more complex nonlinear flow problems,
such as the incompressible Navier--Stokes equations discussed next.

\subsection{Incompressible Navier-Stokes Equation: Flow around a cylinder}

This benchmark is used to evaluate the proposed hyper-reduction strategy in a
realistic CFD setting involving pressure--velocity coupling,
and unsteady flow dynamics. Unlike the previous test cases, this problem
combines strong nonlinearity with a coupled saddle-point structure, making it a
stringent test for hyper-reduced models, particularly when used in conjunction
with segregated solvers.

We consider the unsteady incompressible Navier--Stokes equations
\autoref{eq:navstokes} on a two-dimensional domain representing the flow around
a circular cylinder. The governing equations are discretized using the finite
volume method and solved with a segregated pressure--velocity coupling strategy
based on the SIMPLE algorithm. The computational domain is defined as
$\Omega := [-4D,30D] \times [-6D,6D] \setminus B_D(0,0)$, where $D=1\,\mathrm{m}$
denotes the cylinder diameter.

\begin{equation}
\label{eq:navstokes}
\begin{cases}
\frac{\partial\bm{u}}{\partial t}+ \bm{\nabla} \cdot (\bm{u} \otimes \bm{u}) - \bm{\nabla} \cdot \nu \left(\bm{\nabla}\bm{u}+\left(\bm{\nabla}\bm{u}\right)^T\right)=-\bm{\nabla}p &\mbox{ in } \Omega \times [0,T],\\
\bm{\nabla} \cdot \bm{u}=\bm{0} &\mbox{ in } \Omega \times [0,T],\\
\bm{u} (t,\bm{x} ) = \bm{f}(\bm{x}) &\mbox{ on } \Gamma_\textrm{inlet} \times [0,T],\\
\bm{u} (t,\bm{x} ) = \bm{0} &\mbox{ on } \Gamma_{0} \times [0,T],\\ 
(\nu\bm{\nabla} \bm{u} - p\bm{I})\bm{n} = \bm{0} &\mbox{ on } \Gamma_\textrm{outlet} \times [0,T],\\ 
\bm{u}(0,\bm{x})=\bm{R}(\bm{x}) &\mbox{ in } (\Omega,0),\\            
\end{cases}
\end{equation}
Here,  $t$ is the time, $\bm{x}$ is the spatial variable vector and $\Gamma = \Gamma_\textrm{inlet} \cup \Gamma_0 \cup \Gamma_\textrm{outlet}$ is the boundary of the fluid domain $\Omega$. 
The three parts that form the boundary are called $\Gamma_\textrm{inlet}$, $\Gamma_\textrm{outlet}$, and $\Gamma_0$, they correspond to the inlet boundary, the outlet boundary, and the physical walls, respectively. The fluid kinematic viscosity is denoted by $\nu$ and is constant across the spatial domain. The function $\bm{f}$ includes the boundary conditions for the non-homogeneous boundary. The initial velocity field is given by the function $\bm{R}(x)$. The normal unit vector is denoted by $\bm{n}$. 

\begin{figure}
\centering
\includegraphics[width=\textwidth]{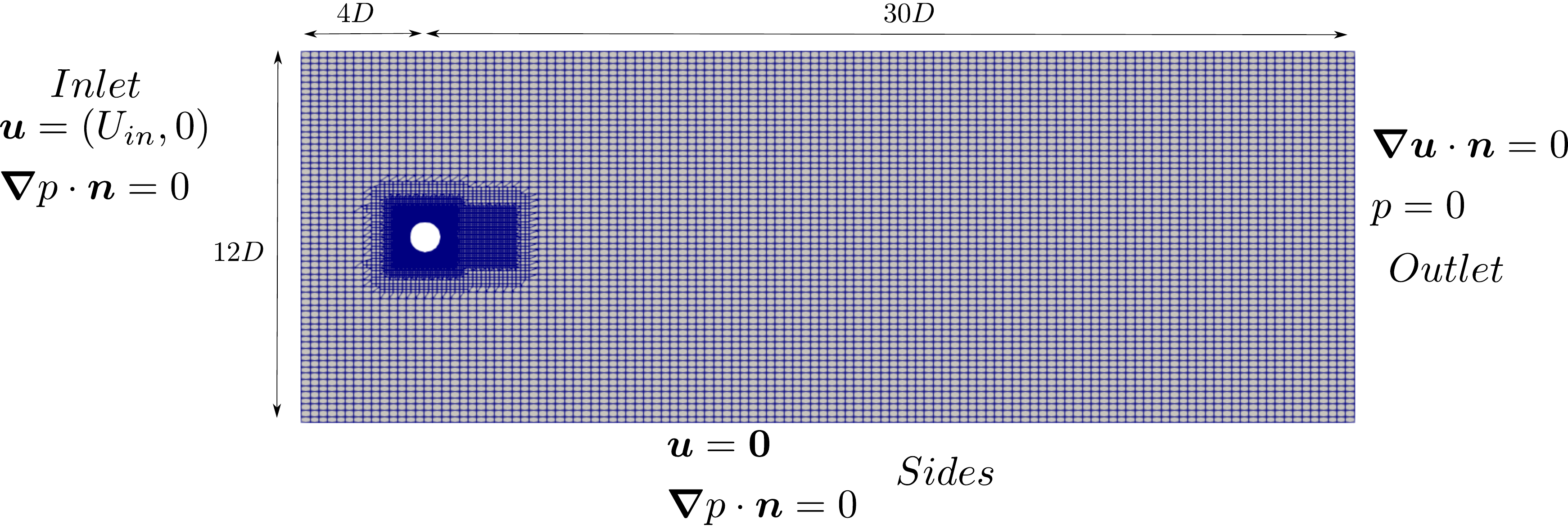} 
\caption{The computational grid of the problem of the flow around a circular cylinder.} 
\label{fig:comp_domain_cy} 
\end{figure}

The mesh, shown in \autoref{fig:comp_domain_cy}, consists of 11\,644 control
volumes and 24\,440 points. A uniform inlet velocity
$\bm{U}_{\infty} = (1.0,0)\,\mathrm{m/s}$ is prescribed, corresponding to a
Reynolds number of $Re=200$, while no-slip conditions are enforced on the
cylinder surface and the channel walls. Zero-gradient conditions are applied at
the outlet.

\begin{figure}
\centering
	\begin{tabular}{cc}
\includegraphics[width=0.75\textwidth]{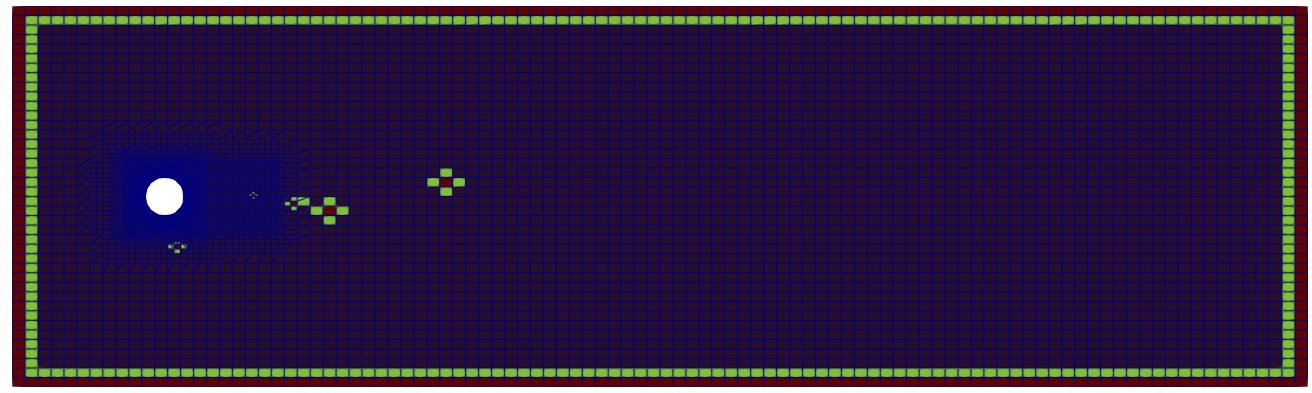}
        &
\includegraphics[width=0.22\textwidth]{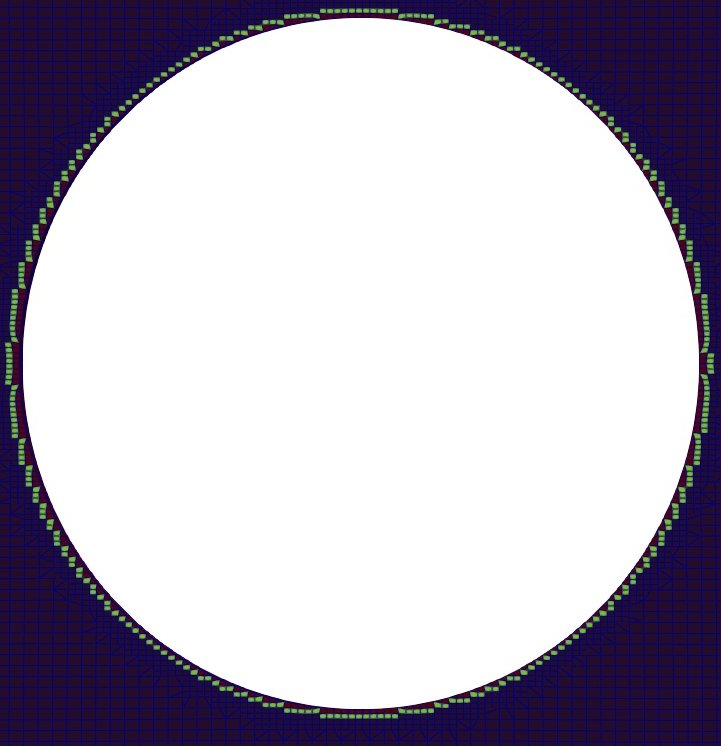}
\\                                                     
(a) & (b)
	\end{tabular}
	\caption{(a) DEIM and obligatory (boundary) cells on the domain, (b) zoom around the obligatory cells on the cylinder. The magic cells are in red, and the yellow cells are the associated neighbours. The discarded nodes in the evolution of the dynamics are in blue. The stencil is made of 1 layer of cells}	
   \label{fig:optim_points_around_cylinder}
\end{figure}

Hyper-reduction is applied independently to the momentum  and pressure equations,
reflecting the segregated nature of the solver. Separate reduced spaces are
constructed for the velocity and pressure fields, with distinct sets of DEIM
sampling points. \autoref{fig:optim_points_around_cylinder} illustrates the
selected interpolation cells together with the associated stencil neighbourhoods.
In addition to DEIM-selected points, obligatory cells on the cylinder boundary
are retained to ensure accurate enforcement of boundary conditions and numerical
stability.

Table~\ref{tab:OfflineNSEs} reports the offline and online computational times for
different combinations of velocity and pressure modes. The observed speed-up
factors, approximately $1.8\times$, are smaller than those obtained for the
linear transport problem. This reduction is expected, as the Navier--Stokes
equations require repeated nonlinear residual evaluations and pressure
corrections at each time step, which limit the achievable acceleration in
segregated formulations.

Despite the modest speed-up, the hyper-reduced model accurately reproduces the
full-order flow dynamics. \autoref{fig:fields_comparison} shows that both the
velocity and pressure fields are well captured at $t=2\,\mathrm{s}$, with
errors remaining small and spatially localized. Importantly, no spurious
pressure oscillations or loss of stability are observed, indicating that the
hyper-reduction strategy is compatible with the pressure--velocity coupling
mechanism of the segregated solver.

This test case demonstrates that the proposed approach can be applied to
pressure-coupled, nonlinear flow problems to
existing finite-volume solvers. While the attainable speed-up is inherently
limited by the structure of the Navier--Stokes equations and the segregated
solution strategy, the method delivers stable and accurate reduced-order
solutions with mesh-independent online complexity, making it suitable for
large-scale industrial CFD applications.

\begin{table}
\centering
\begin{tabular}{|p{2cm}||p{5cm}||p{3cm}||p{2cm}|}
	\hline
	\textbf{Stages} & \textit{\# of modes} & Time [s] & {\RA Speed-up $[T_{Off}/ T_{On}]$} \\
 \hline
	\textbf{Offline} &  - & 166.742 & -\\
 \hline
    \multirow{3}{*}{\textbf{Online}} & $N_{\boldsymbol{u}}=5, N_p=4$ & 88.69 & 1.88 \\
       \cline{2-4}
       &  $N_{\boldsymbol{u}}=3, N_p=5$ & 89.61& 1.86 \\
       \hline
\end{tabular}
\caption{Speed-up, offline and online times comparison varying the number of modes}
\label{tab:OfflineNSEs}
\end{table}

\begingroup
\setlength{\tabcolsep}{0.05pt} 
\begin{figure}
\centering
\begin{tabular}{ccc}
\includegraphics[width=0.33\linewidth]{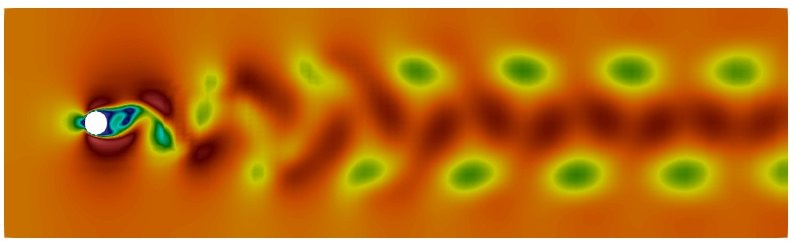} &  \includegraphics[width=0.33\linewidth]{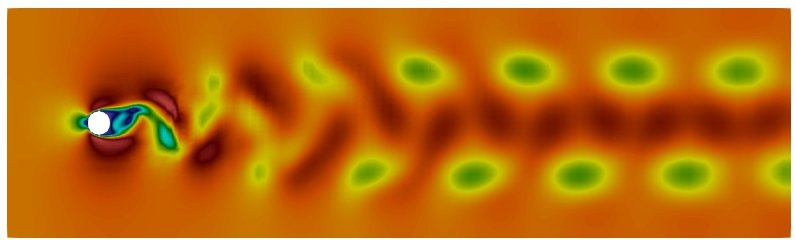} &
\includegraphics[width=0.33\linewidth]{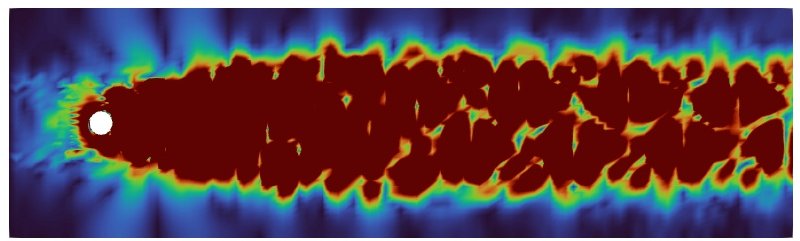}\\
\includegraphics[width=0.23\linewidth]{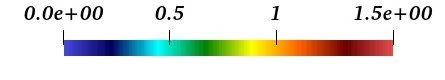}  &  \includegraphics[width=0.23\linewidth]{images/U_scale.jpeg}  & \includegraphics[width=0.23\linewidth]{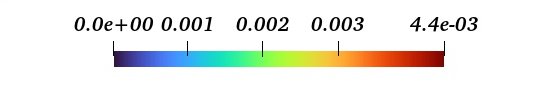} \\
\includegraphics[width=0.33\linewidth]{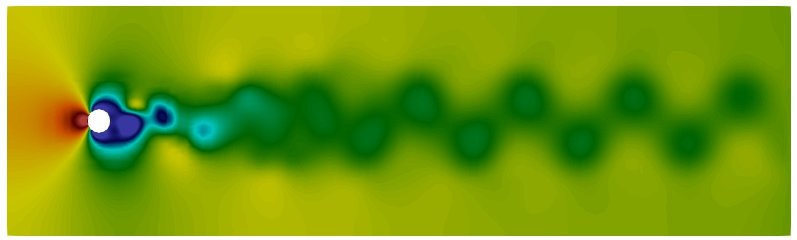} &  \includegraphics[width=0.33\linewidth]{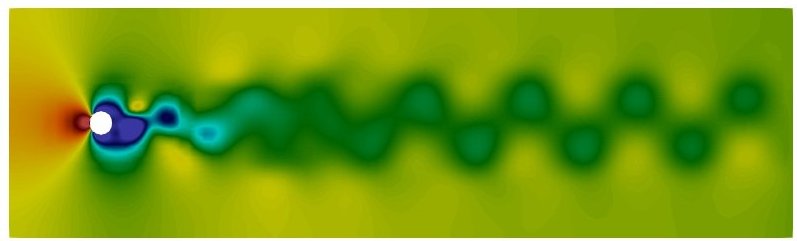} &
\includegraphics[width=0.33\linewidth]{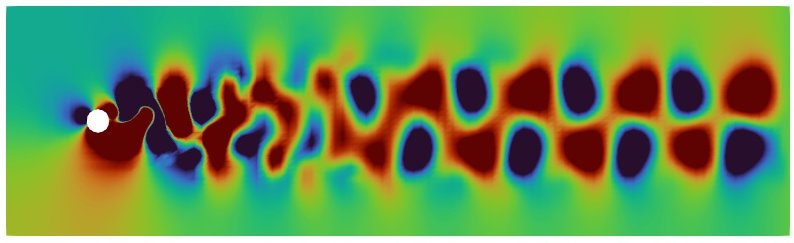}\\
\includegraphics[width=0.23\linewidth]{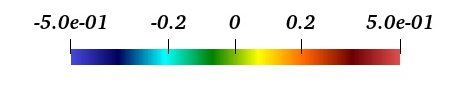}  &  \includegraphics[width=0.23\linewidth]{images/p_scale.jpeg}  & \includegraphics[width=0.23\linewidth]{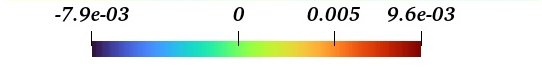}  \\
(a) & (b) & (c)
\end{tabular}
\caption{Comparison of the velocity (first row) and pressure (second row) fields. 
The first column shows the original solutions, the second column shows the hyper-reduced  solutions obtained using the first 3 DEIM points for $\bm{U}$ and 5 DEIM points 
for $p$, respectively. The third column depicts the associated errors.  All snapshots are taken at $t = 2\,\mathrm{s}$}	
\label{fig:fields_comparison}
\end{figure}
\endgroup

\subsection{Discussion}

Across all test cases, the proposed hyper-reduction strategy achieves
mesh-independent online costs while preserving the structure of segregated
finite-volume solvers. Linear problems benefit from significant speed-ups,
whereas nonlinear and pressure-coupled systems exhibit more modest gains,
reflecting the intrinsic cost of nonlinear residual evaluations. Importantly,
the method remains stable and non-intrusive, even for turbulent flows with
complex boundary conditions, which are typically challenging for standard
hyper-reduction techniques.

\section{Conclusions and perspectives}
\label{conclusion}

This work introduced a hyper-reduced order modeling strategy tailored to
segregated finite-volume solvers and geometrically parametrized problems.
The method follows an approximate-then-project paradigm: the full-order
operators are assembled once and subsequently evaluated online using a
small, carefully selected subset of mesh cells. This design removes the
dependence of the online computational cost on the full mesh size while
preserving the structure of the underlying solver.

The proposed approach was validated on three representative problems: a
linear transport equation, a nonlinear Burgers equation, and the
incompressible Navier--Stokes equations for flow around a cylinder. In all
cases, the hyper-reduced models accurately reproduced the full-order
solutions. For linear problems, substantial speed-ups were achieved, while
for nonlinear and pressure-coupled systems the acceleration was more
moderate but remained consistent with the intrinsic cost of nonlinear
residual evaluations. Importantly, stability and accuracy were maintained
across all test cases.

A key feature of the method is its compatibility with segregated finite-volume
solvers and its reliance on local stencil information only. This makes the
approach non-intrusive and directly applicable to existing industrial CFD
codes, even in the presence of geometric variations and mesh motion.
Since only a sparse subset of cells is processed during the online phase,
the method is naturally parallel and scalable, making it suitable for
large-scale simulations on modern high-performance computing architectures.

Future work will focus on adaptive sampling strategies, extension to
strongly coupled multiphysics problems, and further reduction of nonlinear
operator costs. These directions aim to enhance the efficiency of the
approach while preserving its solver-aware and geometry-compatible
characteristics.

\section*{Acknowledgements}

Gianluigi Rozza acknowledges funding from the European Union Horizon 2020 Program in the framework of European Research Council Executive Agency: H2020 ERC CoG 2015 AROMA-CFD project 681447 "Advanced Reduced Order Methods with Applications in Computational Fluid Dynamics" P.I. Professor Gianluigi Rozza. Gianluigi Rozza also acknowledges funding from the Italian Ministry of University and Research (MUR) in the form of PRIN ``Numerical Analysis for Full and Reduced Order Methods for Partial Differential Equations'' (NA-FROM-PDEs) project, and by INdAM GNCS.
Andrea Mola acknowledges the financial support under the National Recovery and Resilience Plan (NRRP), Mission 4, Component 2, Investment 1.1, Call for tender No. 1409, funded by the European Union – NextGenerationEU – Project Title ROMEU – CUP D53D2301888001 - Grant Assignment Decree No. 1379 adopted on 01/09/2023 by the Italian Ministry of University and Research (MUR) and the financial support by INdAM GNCS.
Giovanni Stabile acknowledges the financial support under the National Recovery and Resilience Plan (NRRP), Mission 4, Component 2, Investment 1.1, Call for tender No. 1409, funded by the European Union – NextGenerationEU– Project Title ROMEU – CUP P2022FEZS3 - Grant Assignment Decree No. 1379 adopted on 01/09/2023 by the Italian Ministry of University and Research (MUR), and acknowledges the financial support by the European Union (ERC, DANTE, GA-101115741). Views and opinions expressed are however, those of the author(s) only and do not necessarily reflect those of the European Union or the European Research Council Executive Agency. Neither the European Union nor the granting authority can be held responsible for them.
\bibliographystyle{abbrv}
\bibliography{biblio}
\end{document}